\documentclass{amsart}

\usepackage[mathcal]{euscript}
\usepackage{amssymb, amsfonts, xypic}

\newtheorem{theorem}{Theorem}[section]
\newtheorem{lemma}[theorem]{Lemma}

\newtheorem{proposition}[theorem]{Proposition}

\theoremstyle{definition}
\newtheorem{definition}[theorem]{Definition}

\theoremstyle{remark}
\newtheorem{remark}[theorem]{Remark}

\numberwithin{equation}{section}

\newcommand{\C}{\mathcal{C}}

\DeclareMathOperator{\Ho}{Ho}
\DeclareMathOperator{\Hom}{Hom}

\DeclareMathOperator*{\colim}{colim}

\newcommand{\pros}{\textup{pro}}
\newcommand{\pro}{\textup{pro-}}
\newcommand{\Map}{\textup{Map}}

\newcommand{\map}{\rightarrow}

\newcommand{\tl}{\tilde}

\newcommand{\dfn}{\textbf} 
\newcommand{\mdfn}[1]{\dfn{\mathversion{bold}#1}} 

\begin{document}

\title{Strict Model Structures for Pro-Categories}

\author{Daniel C. Isaksen}
\address{Department of Mathematics \\ University of Notre Dame \\
Notre Dame, IN 46556}
\email{isaksen.1@nd.edu}

\subjclass{18G55, 55U35}
\date{August 24, 2001}
\keywords{Closed model structures, pro-homotopy theory, pro-spaces}

\thanks{}

\begin{abstract}
We show that if $\C$ is a proper model category, then the 
pro-category $\pro \C$ has a strict model structure in which the 
weak equivalences are the levelwise weak equivalences.  This is related
to a major result of \cite{EH}.  The strict model structure is the
starting point for many homotopy theories of pro-objects such as those
described in \cite{duality}, \cite{prospace}, and \cite{completion}.
\end{abstract}

\maketitle

\section{Introduction}
\label{sctn:intro}

If $\C$ is a category, then the category $\pro \C$ 
has as objects all
cofiltered diagrams in $\C$ and has morphisms defined by
\[
\Hom_{\pro \C} (X, Y) = \lim_t \colim_s \Hom_{\C} (X_s, Y_t).
\]
Pro-categories have found many uses over the years in fields such 
as algebraic geometry \cite{AM}, shape theory \cite{MS}, 
geometric topology \cite{CJS}, and possibly even applied 
mathematics \cite[Appendix]{CP}.

When working with pro-categories, one would frequently like to have
a homotopy theory of pro-objects.  The first attempts at this
appear in \cite{AM} and \cite{Sullivan} in which pro-objects in
homotopy categories are considered.  The difficulty with this approach
is that the diagrams commute only up to homotopy, and this makes it
virtually impossible to make sense of most of the standard notions
of homotopy theory in this context.

Much better is to first consider actually commuting cofiltered diagrams
(of spaces or simplicial sets or spectra or whatever) and then to define
a notion of weak equivalence between such pro-objects.  This approach
was first taken by \cite{JG} in a restricted context.  

It was
also applied much more generally in \cite{EH}.  The idea is to start
with a model structure ({\em i.e.}, a homotopy theory) on a category $\C$ 
and then to construct a 
{\em strict model structure}
on $\pro \C$ in which the weak equivalences are more or less just the
levelwise weak equivalences.  The resulting homotopy theory 
is precisely suited
to study homotopy limits \cite[Ch.~XI]{BK} of cofiltered diagrams.
The strict model structure is a starting point for several
other model structures such as those described in \cite{duality}, 
\cite{prospace}, and \cite{completion}.  

The strict model structure on $\pro \C$ does not seem
to exist for a completely arbitrary model category $\C$.  A niceness
hypothesis was required in \cite[p.~45]{EH} (which was weakened in
\cite{HH}).
Unfortunately, many important examples of model categories,
such as the usual models for spectra, 
do not satisfy this hypothesis.
The main purpose of this paper is to prove that the strict structure
on $\pro \C$ exists whenever $\C$ is a proper model category.  
Almost all of the most important examples of model categories 
are proper.

Another problem with \cite{EH} is that a non-standard set of axioms
for model structures are used.  From a modern perspective, it is
harder to comprehend the technical details of \cite{EH} than to simply
work out new proofs from scratch.  The secondary goal of this paper
is to give these new more modern proofs.
Oddly, the two-out-of-three axiom is the most
difficult part of the proof; 
in most model structures, it is automatic from the definition
of weak equivalence.  

The last goal of the paper is to consider whether strict model structures
on pro-categories are fibrantly generated.  
It is already known that these model structures are not cofibrantly
generated in general, even when $\C$ is \cite[\S 19]{prospace}.
We produce reasonable
collections of generating fibrations and generating acyclic fibrations
that have cosmall codomains, 
but these collections are not sets.  
In fact, the strict model structure for pro-simplicial sets is not
fibrantly generated.  We show that if this strict structure were
fibrantly generated, then in the category of simplicial sets 
there would exist a set of fibrations that detect acyclic cofibrations.

The paper is organized as follows.  First we introduce the language
of pro-categories and give some background results.  Then we
define the strict weak equivalences and prove that they
satisfy the two-out-of-three axiom when $\C$ is proper.  
Next we prove that the strict model structure exists when $\C$
is proper.  Finally, we consider whether the strict model structure
is fibrantly generated.

We assume familiarity with model categories.  The original
reference is \cite{Quillen}, but we follow the notation and terminology
of \cite{PH} as closely as possible.
Other references include \cite{Dwyer-Spalinski} and \cite{Hovey}.

\section{Preliminaries on Pro-Categories} \label{sctn:prelim-pro}

We begin with a review of the necessary background on pro-categories.
This material can be found in 
\cite{SGA}, 
\cite{AM}, \cite{Duskin}, \cite{EH}, and \cite{prolimits}.

\subsection{Pro-Categories}

\begin{definition} 
\label{defn:pro}
For a category $\C$, the category \mdfn{$\pro \C$} has objects all 
cofiltering
diagrams in $\C$, and 
$$\Hom_{\pro \C}(X,Y) = \lim_s \colim_t \Hom_{\C}
     (X_t, Y_s).$$
Composition is defined in the natural way.
\end{definition}

A category $I$ is \dfn{cofiltering} if the following conditions hold:
it is non-empty and small;
for every pair of objects $s$ and $t$ in $I$,
there exists an object $u$ together with maps $u \map s$ and
$u \map t$; and for every pair of morphisms $f$ and $g$ with the
same source and target, there exists a morphism $h$ such that $fh$ equals
$gh$.
Recall that a category is \dfn{small} if it has only a set 
of objects and a set of morphisms.
A diagram is said to be \dfn{cofiltering} if its
indexing category is so.
Beware that some material on pro-categories, 
such as \cite{AM} and \cite{Meyer},
consider cofiltering categories that are not small.
All of our pro-objects will be indexed by small categories.

Objects of $\pro \C$ are functors from
cofiltering categories to $\C$.  
We use both set theoretic and categorical
language to discuss indexing categories;
hence ``$t \geq s$'' and ``$t \map s$'' mean the same thing 
when the indexing category is actually a directed set.

The word {\em pro-object} refers to objects of pro-categories.
A \dfn{constant} pro-object is one indexed
by the category with one object and one (identity) map.
Let $c: \C \map \pro \C$ be the functor taking an object $X$ to the 
constant pro-object with value $X$.
Note that this functor makes $\C$ a full subcategory of
$\pro \C$.  The limit functor $\lim: \pro \C \map \C$
is the right adjoint of $c$.  
To avoid confusion, we write $\lim^{\pros}$ for limits computed 
within the category $\pro \C$.

Let $Y: I \map \C$ and $X: J \map \C$ be arbitrary pro-objects.
We say that $X$ is 
\dfn{cofinal} in $Y$ if there is a cofinal functor
$F: J \map I$ such that
$X$ is equal to the composite $Y F$.
This means that for every $s$ in $I$,
the overcategory $F \downarrow s$ is cofiltered.
In the case when $F$ is an inclusion
of directed sets, $F$ is cofinal if and only if
for every $s$ in $I$ there
exists $t$ in $J$ such that $t \geq s$.
The importance of this definition is that $X$ is isomorphic to $Y$ in
$\pro \C$.

A \dfn{level representation} of a map
$f:X \map Y$ is:
a cofiltered index category $I$;
pro-objects $\tl{X}$ and $\tl{Y}$ indexed by $I$ and pro-isomorphisms
$X \map \tl{X}$ and $Y \map \tl{Y}$;
and a collection of
maps $f_{s}:\tl{X}_s \map \tl{Y}_s$ for all $s$ in $I$
such that for all $t \map s$ in $I$, there is a commutative diagram
\[
\xymatrix{
\tl{X}_{t} \ar[d] \ar[r] & \tl{Y}_{t} \ar[d] \\
\tl{X}_{s} \ar[r]        & \tl{Y}_{s}       }
\]
and such that the maps $f_{s}$ represent a map $\tl{f}: \tl{X} \map \tl{Y}$
belonging to a commutative square
\[
\xymatrix{
X \ar[d] \ar[r]^{f} & Y \ar[d] \\
\tl{X} \ar[r]_{\tl{f}} & \tl{Y}     }
\]
in $\pro \C$.
That is, a level representation is just a natural transformation 
such that the maps $f_{s}$ represent the element $f$ of
\[
\lim_{s} \colim_{t} \Hom_{\C} (X_{t}, Y_{s}) \cong
\lim_{s} \colim_{t} \Hom_{\C} (\tl{X}_{t}, \tl{Y}_{s}).
\]
Every map has a level representation 
\cite[App.~3.2]{AM} \cite{Meyer}.

More generally, suppose given any diagram $A \map \pro \C: a \mapsto X^{a}$.
A \dfn{level representation}
of $X$ is: a 
cofiltered index category $I$; a 
functor $\tl{X}: A \times I \map \C: (a, s) \mapsto \tl{X}^{a}_{s}$;
and pro-isomorphisms $X^{a} \map \tl{X}^{a}$ such that 
for every map $\phi:a \map b$ in $A$, 
$\tl{X}^{\phi}$ is a level representation for $X^{\phi}$.
In other words, $\tl{X}$ is a uniform level representation for all the maps
in the diagram $X$.

Not every diagram of pro-objects has a level representation.  However,
finite diagrams without loops do have level representations.  This makes
computations of limits and colimits of such diagrams in $\pro \C$ 
relatively straightforward.  To compute this limit or colimit,
just take the levelwise limit or colimit
of the level representation \cite[App.~4.2]{AM}.

A pro-object $X$ satisfies a certain property \dfn{levelwise} if
each $X_{s}$ satisfies that property, and $X$ satisfies this property
\dfn{essentially levelwise} if it is isomorphic to another pro-object
satisfying this property levelwise.
Similarly, a level representation $X \map Y$ 
satisfies a certain property \dfn{levelwise}
if each $X_{s} \map Y_{s}$ has this property.
A map of pro-objects satisfies this property \dfn{essentially levelwise}
if it has a
level representation satisfying this property levelwise.
The following 
surprisingly general and very useful proposition about retracts of 
essentially levelwise maps is proved in \cite[Thm.~5.5]{prolimits}.

\begin{proposition} 
\label{prop:retract}
Let $C$ be any class of maps in a category $\C$.  Then retracts preserve
the class of maps in $\pro \C$ that belong to $C$ essentially levelwise.
\end{proposition}

\subsection{Cofiniteness}

A directed set $(I,\leq)$ is \dfn{cofinite}
if for every $t$, the set of elements $s$ of $I$ such that $s \leq t$
is finite.  
A pro-object or level representation is \dfn{cofinite directed} 
if it is indexed by a cofinite directed set.

For every cofiltered category $I$, there exists a cofinite directed
set $J$ and a cofinal functor $J \map I$
\cite[Th.~2.1.6]{EH} 
(or \cite[Expos\'e~1, 8.1.6]{SGA}).  
Therefore, every
pro-object is isomorphic to a cofinite directed pro-object.
Similarly, every map has
a cofinite directed level representation.
Thus, it is possible to restrict the definition of a pro-object to
only consider cofinite directed sets as index categories, but 
we find this unnatural for general definitions and constructions.
On the other hand, we find it much easier to work with
cofinite directed pro-objects in practice.  Thus, most of our results
start by assuming without loss of generality 
that a pro-object is indexed by a cofinite directed set.
Cofiniteness is critical because many arguments and constructions proceed
inductively.

\begin{definition} \label{defn:matching-map}
Let $f: X \map Y$ be a cofinite directed level representation of
a map in a pro-category.  For every index $t$, the 
\mdfn{relative matching map $M_{t} f$} is the map 
\[
X_{t} \map \lim_{s<t} X_{s} \times_{\lim_{s<t} Y_{s}} Y_{t}.
\]
\end{definition}

The terminology is motivated by the fact that these maps
appear in Reedy model structures \cite[Defn.~16.3.2]{PH}.
The similarity is not coincidental.  The
strict model structure is
closely linked to the Reedy model structures for each fixed cofinite
directed index category \cite[\S 3.2]{EH}.

\section{Strict Weak Equivalences} \label{sctn:strict}

We now study strict weak equivalences for pro-categories
as originally described in \cite{EH}.
The niceness hypothesis of \cite[p. 45]{EH} is not satisfied by
many categories of interest.  These include many of the standard
models for spectra, such as Bousfield-Friedlander spectra \cite{BF},
symmetric spectra \cite{HSS}, or $S$-modules \cite{EKMM}.
We shall study strict weak equivalences in $\pro \C$ whenever
$\C$ is a proper model category.  

\begin{definition}
\label{defn:strict-we}
The \dfn{strict weak equivalences} 
of $\pro \C$ are the essentially
levelwise weak equivalences.
\end{definition}

Actually, a more complicated definition appears in \cite[\S~3.3]{EH},
but we shall prove in Proposition \ref{prop:strict-equivalent} 
the equivalence of that definition and ours.

It is not obvious from the definition that the strict weak equivalences
satisfy the two-out-of-three axiom.  
The next few lemmas prove this axiom.
These proofs are the technical
heart of the paper.  They are the reason that we must assume that $\C$
is proper.  The basic complication is that given a diagram
of strict weak equivalences, it is not necessarily possible to 
find a level representation for the diagram such that the 
level representations of all maps are levelwise weak equivalences.
Independently, each map has a level representation that is a levelwise
weak equivalence, but the reindexing for these level representations
may be different.

\begin{lemma} \label{lem:pro-factor}
Let $\C$ be a model category, and let $f:X \map Y$ be a level representation
of an isomorphism in $\pro \C$.  After reindexing along a cofinal
functor, $f$ can be
factored into a levelwise cofibration that is also
a pro-isomorphism followed by a
levelwise fibration that is also a pro-isomorphism.
\end{lemma}

\begin{remark}
\label{rem:pro-factor}
The model structure on $\C$ is not really necessary.  We just need a category
$\C$ and two classes of morphisms $C$ and $F$ such that each morphism of 
$\C$ can be factored (not necessarily functorially) into an element
of $C$ followed by an element of $F$.
\end{remark}

\begin{proof}
We may assume that $f$ is
indexed by a directed set $I$.
Since $f$ is an isomorphism, for every $s$
in $I$, there exists $t > s$ and a map $h_{ts}:Y_t \map X_s$ belonging
to a commutative diagram
\[
\xymatrix{
X_t \ar[r] \ar[d] & Y_t \ar[d] \ar[dl] \\
X_s \ar[r] & Y_s.                        }
\]
In effect, the maps $h_{ts}$ represent the inverse of $f$.
By restricting to cofinal subsets, we may assume that such a diagram
exists for every $t > s$.  We choose a map $h_{ts}: Y_t \map X_s$
for each pair $t > s$.  

Factor each map $f_s: X_s \map Y_s$ into a cofibration
$X_s \map Z_s$ followed by a fibration $Z_s \map Y_s$.
We shall define structure maps making $Z$ into
a pro-object.

Define a category $J$ whose objects are the elements of $I$ and whose
morphisms $t \map s$ are finite chains $t = u_0 > u_1 > \cdots > u_n = s$.
Composition is defined by concatenation of chains.  Note that
$J$ is not a directed set; it is not even cofiltered.

For every morphism $t = u_0 > u_1 > \cdots > u_n = s$ in $J$, define
a map $Z_t \map Z_s$ by the composition
\[
\xymatrix@1{
Z_t \ar[r] & Y_t = Y_{u_0} \ar[r]^{h_{u_{0}u_{1}}} & X_{u_1} \ar[r] &
   X_{u_2} \ar[r] & \cdots \ar[r] & X_{u_n} = X_s \ar[r] & Z_s.   }
\]
A diagram chase shows that this makes $Z$ into a diagram indexed by $J$.

There may
be more than one map from a given $Z_t$ to a given $Z_s$, but another
diagram chase shows that they become equal after composition with
some map $Z_u \map Z_t$.
Consider the category $K$ defined to be a quotient of $J$ as follows.
The objects of $K$ are the same as the objects of $J$, and two morphisms
from $t$ to $s$ in $J$ are identified in $K$ if the corresponding
maps from $Z_t$ to $Z_s$ are equal in $\C$.
Now $K$ is a cofiltered category, and we may consider it as the indexing
category of $Z$.

The projection functor $K \map I$ is cofinal, so
we may reindex $X$ and $Y$ along this functor.
More diagram chases show that the maps $X_s \map Z_s$ and $Z_s \map Y_s$
assemble into level representations $X \map Z$ and $Z \map Y$.
It remains only to show that these maps are pro-isomorphisms.
This follows from the commutative diagrams
\[
\xymatrix{
X_{t} \ar[r] \ar[d] & Z_{t} \ar[d] \ar[dl] &  
    Z_{t} \ar[r] \ar[d] & Y_{t} \ar[d] \ar[dl] \\
X_{s} \ar[r] & Z_{s} & 
    Z_{s} \ar[r] & Y_{s}          }
\]
for every pair $t > s$.
The diagonal maps in the above diagrams are the compositions
$Z_t \map Y_t \map X_s$ and 
$Y_t \map X_s \map Z_s$ respectively.
\end{proof}

\begin{remark}
\label{rem:pro-factor2}
There is a more obvious argument where a level representation $X \map Y$
is functorially factored into a levelwise cofibration $X \map Z$ followed
by a levelwise fibration $Z \map Y$.  This does not give the desired 
factorization; the structure maps for $Z$ are wrong. 
\end{remark}

The following lemma appears in \cite[Prop.~10.4]{prospace}, 
but the previous lemma
makes the technical details of that proof clearer.

\begin{lemma} \label{lem:composition}
If $\C$ is a proper model category, then the strict weak 
equivalences of $\pro \C$ are closed under composition.
\end{lemma}

\begin{proof}
It suffices to 
assume that there is a cofinite directed level representation for the diagram
$$\xymatrix{X \ar[r]^f & Y & Z \ar[l]_h \ar[r]^g & W}$$
in which $f$ and $g$ are levelwise weak equivalences while
$h$ is a pro-isomorphism (but not a levelwise isomorphism).
We must construct a levelwise weak equivalence isomorphic to the 
composition $gh^{-1}f$.

By Lemma \ref{lem:pro-factor}, after reindexing we can 
factor $h: Z \map Y$ into a levelwise
cofibration $Z \map A$ followed by a levelwise fibration $A \map Y$
such that 
\[
X \map Y \leftarrow A \leftarrow Z \map W
\]
is a level representation in which the first and fourth maps
are levelwise weak equivalences, 
and the second and third are
pro-isomorphisms.

Let $B$ be the pullback $X \times_Y A$, 
and let $C$ be the pushout $A \amalg_Z W$.
Since pushouts and pullbacks can be constructed levelwise, the
maps $B \map A$ and $A \map C$ are levelwise weak equivalences.  Here
we use that the model structure is left and right proper.
Moreover, the maps $B \map X$ and $W \map C$ are pro-isomorphisms
since base and cobase changes preserve isomorphisms.
Hence the composition $B \map C$ is the desired levelwise weak equivalence.
\end{proof}

\begin{lemma} \label{lem:2/3}
Suppose that $\C$ is a proper model category, and 
let $f$ and $g$ be two composable maps in pro-$\C$.
If $g$ and $gf$ are strict weak equivalences,
then $f$ is a strict weak equivalence.
If $f$ and $gf$ are strict weak equivalences,
then $g$ is a strict weak equivalence.
\end{lemma}

\begin{proof}
We first prove the first claim.
We may consider a cofinite directed level representation
\[
\xymatrix{
X \ar[r]^f \ar[d] & Y \ar[d]^g \\
W \ar[r] & Z}
\]
where the vertical maps are levelwise
weak equivalences and the bottom horizontal map is a pro-isomorphism
(but not a levelwise isomorphism).
We want to show that the top horizontal map is an essentially
levelwise  weak equivalence.

By Lemma \ref{lem:pro-factor}, after reindexing
there exists a level representation
\[
\xymatrix{
X \ar[r] \ar[d] & B \ar[r] \ar[d] & Y \ar[d] \\
W \ar[r] & A \ar[r] & Z}
\]
such that $B$ is the levelwise pullback $A \times_Z Y$;
the first and third vertical maps are levelwise weak equivalences;
the map $W \map A$ is a levelwise {\em acyclic}
cofibration and a pro-isomorphism;
and the map $A \map Z$ is a levelwise fibration and a pro-isomorphism.
Because of right properness, the map $B \map A$ is also a levelwise weak
equivalence.
By the two-out-of-three axiom in $\C$, the induced map
$X \map B$ is a levelwise weak equivalence.

On the other hand, the map $B \map Y$ is an isomorphism because
base changes preserve isomorphisms.
Hence $X \map B$ is isomorphic to $f$.

The proof of the second claim is similar.  We start with a
cofinite directed level representation
\[
\xymatrix{
X \ar[r] \ar[d]_f & W \ar[d] \\
Y \ar[r]_g & Z}
\]
where the vertical maps are levelwise
weak equivalences and the top horizontal map is a pro-isomorphism.
We want to show that the bottom horizontal map is an essentially
levelwise weak equivalence.
We produce a level representation
\[
\xymatrix{
X \ar[r] \ar[d] & A \ar[r] \ar[d] & W \ar[d] \\
Y \ar[r] & B \ar[r] & Z}
\]
such that $B$ is the levelwise pushout $A \amalg_X Y$;
the first and third vertical maps are levelwise weak equivalences;
the map $A \map W$ is a levelwise {\em acyclic}
fibration and a pro-isomorphism;
and the map $X \map A$ is a levelwise cofibration and a pro-isomorphism.
The map $B \map Z$ is the desired level representation.
\end{proof}

The above lemmas imply that
our definition of strict weak equivalences agrees with
the definition of \cite[\S 3.3]{EH}.

\begin{proposition} \label{prop:strict-equivalent}
When $\C$ is a proper model category,
the strict weak equivalences of Definition \ref{defn:strict-we} agree
with the weak equivalences of \cite[\S 3.3]{EH}.
\end{proposition}

\begin{proof}
The weak equivalences of \cite[\S 3.3]{EH} are by definition compositions of 
essentially levelwise weak equivalences.  By Lemma~\ref{lem:composition},
these compositions are again essentially levelwise weak equivalences.
On the other hand, every levelwise weak equivalence can be factored into
a levelwise acyclic cofibration followed by a levelwise acyclic fibration.
Therefore, every levelwise weak equivalence is a weak equivalence in
the sense of \cite[\S 3.3]{EH}.
\end{proof}

The preceding proposition is closely related to the main result of 
\cite{Porter}.  However, we make a useful observation missed there.  
Namely, it is not necessary to saturate the essentially levelwise
weak equivalences; they are already saturated when $\C$ is proper.

\section{Strict Model Structures}

Beginning with a proper model structure on a category $\C$, we now
establish a model structure on the category $\pro \C$.

\begin{definition}
\label{defn:cofib}
The \dfn{strict cofibrations} of $\pro \C$ are the
essentially levelwise cofibrations.  
\end{definition}

\begin{definition}
\label{defn:fib}
A map in $\pro \C$ is a \dfn{special fibration} if it has
a cofinite directed level representation $p$ for which every relative matching
map $M_s p$ is a fibration.
A map in $\pro \C$ is a \dfn{fibration} if it
is a retract of a special fibration.
\end{definition}

In order to help us understand these definitions, we need some
auxiliary notions.

\begin{definition}
\label{defn:spec-a-fib}
A map in $\pro \C$ is a \dfn{special acyclic fibration} 
if it has
a cofinite directed level representation $p$ for which every relative matching
map $M_s p$ is an acyclic fibration.
\end{definition}

\begin{remark}
\label{rem:spec-a-fib}
Every special acyclic fibration is a special fibration, so it is
also a strict fibration.  We shall see below in Proposition 
\ref{prop:a-fib} that
the class of strict acyclic fibrations is equal to the class of
retracts of special acyclic fibrations.
\end{remark}

\begin{lemma}
\label{lem:spec-a-fib}
Special acyclic fibrations are essentially levelwise acyclic fibrations.
In particular, they are strict acyclic fibrations.
\end{lemma}

\begin{proof}
Special acyclic fibrations are special fibrations, so they are
strict fibrations by definition.  This means that the second statement
follows from the first.

Suppose given a cofinite directed
level representation $p: X \map Y$ for which each
$M_s p$ is an acyclic fibration.  The map $p_s: X_s \map Y_s$
factors as 
\[
\xymatrix@1{
X_s \ar[r]^-{M_s p} & Y_s \times_{\lim_{t<s} Y_t} \lim_{t<s} X_t
\ar[r]^-{q_s} & Y_s.                                        }
\]
Since compositions and base changes preserve acyclic fibrations, 
it suffices to show inductively in $s$ that 
$\lim_{t<s} p_t: \lim_{t<s} X_t \map \lim_{t<s} Y_t$
is an acyclic fibration for every $s$.
Using that $M_t p$ is an acyclic fibration for $t < s$, 
an induction in $t$ shows that
$\lim_{t<s} p_t$ is also an acyclic fibration.
\end{proof}

\begin{lemma}
\label{lem:level-fib}
Strict fibrations are essentially levelwise fibrations.
\end{lemma}

\begin{proof}
The proof of Lemma \ref{lem:spec-a-fib} shows that special
fibrations are essentially levelwise fibrations.
Retracts of special fibrations are also essentially levelwise fibrations
by Proposition \ref{prop:retract}.
\end{proof}

\begin{lemma}
\label{lem:factor1}
Every map $f:X \map Y$ in $\pro \C$ factors as a strict cofibration 
$i:X \map Z$ followed by
a special acyclic fibration $p:Z \map Y$.
\end{lemma}

\begin{proof}
We may suppose that $f$ is a level representation indexed by a cofinite
directed set.  Suppose for induction
that the maps $i_t: X_t \map Z_t$ and $p_t:Z_t \map Y_t$ have already
been defined for $t < s$.  
Consider the map 
\[
X_s \map Y_s \times_{\lim_{t<s} Y_t} \lim_{t<s} Z_t.
\]
Factor it into a cofibration $i_s: X_s \map Z_s$ followed by an
acyclic fibration
\[
p_s: Z_s \map Y_s \times_{\lim_{t<s} Y_t} \lim_{t<s} Z_t.
\]
This extends the factorization to level $s$.
\end{proof}

\begin{lemma}
\label{lem:factor2}
Every map $f:X \map Y$ in $\pro \C$ factors as an essentially levelwise
acyclic cofibration
$i:X \map Z$ followed by
a special fibration $p:Z \map Y$.
\end{lemma}

\begin{proof}
The proof is identical to the proof of Lemma \ref{lem:factor1}, except
that we factor the map
\[
X_s \map Y_s \times_{\lim_{t<s} Y_t} \lim_{t<s} Z_t
\]
into an acyclic cofibration followed by a fibration.
\end{proof}

\begin{remark}
\label{rem:a-cofib}
The previous lemma uses the class of
essentially levelwise acyclic cofibrations.
It does
not follow immediately from the definitions that this class is equal
to the class of strict acyclic cofibrations.  If a map is isomorphic
to a levelwise weak equivalence and to a levelwise cofibration, it may
be necessary to use two different index categories to obtain the
two level representations.
In fact, the essentially levelwise acyclic cofibrations 
are the same as the strict acyclic cofibrations, 
as we shall show in Proposition \ref{prop:a-cofib}.
\end{remark}

\begin{remark}
\label{rem:functorial}
The previous two lemmas do not give functorial factorizations.
The problem is that the reindexing into cofinite directed level
representations is not functorial.  This is the reason that
strict model structures do not have functorial factorizations.
\end{remark}

Next we show that the classes of strict cofibrations and
retracts of special acyclic fibrations determine each other by
lifting properties.  Similarly,
the classes of essentially levelwise acyclic cofibrations and
strict fibrations determine each other.

\begin{lemma}
\label{lem:lift1}
A map is a strict cofibration if and only if it has the left lifting
property with respect to all retracts of special acyclic fibrations.
Also, a map has the right lifting property with respect
to all strict cofibrations if and only if it is a retract of a
special acyclic fibration.
\end{lemma}

\begin{proof}
Suppose given a square 
\[
\xymatrix{
A \ar[r] \ar[d]_{i} & X \ar[d]^p \\
B \ar[r] & Y                      }
\]
in which $i$ is a strict cofibration and $p$ is a special acyclic
fibration.
By reindexing,
we may assume that the
diagram is a cofinite directed level representation such that each
$i_s: A_s \map B_s$ is a cofibration and such that each map
$M_s p$ is an acyclic fibration.

Assume that for all $t <s$, we have already
constructed maps $B_{a(t)} \map X_t$ belonging to a commuting diagram
\[
\xymatrix{
A_{a(t)} \ar[r] \ar[d] & X_t \ar[d] \\
B_{a(t)} \ar[r] \ar[ur] & Y_t \times_{\lim_{u<t} Y_u} \lim_{u<t} X_u.     }
\]
Choose $a(s)$ to be
greater than $a(t)$ for all $t <s$ and also greater than $s$; 
this is possible because the indexing
set is cofinite.  Thus we have a diagram
\[
\xymatrix{
A_{a(s)} \ar[r] \ar[d]_{i_{a(t)}} & X_s \ar[d]^{M_s p} \\
B_{a(s)} \ar[r] & Y_s \times_{\lim_{t<s} Y_s} \lim_{t<s} X_s.     }
\]
A lift exists in this diagram since $i_{a(t)}$ is a cofibration and
$M_s p$ is an acyclic fibration.  
By induction on $s$, we construct a lift $B \map X$.

We have shown that strict cofibrations lift with respect to
special acyclic fibrations.
Formally, it follows that strict cofibrations lift
with respect to retracts of special acyclic fibrations.

Now suppose that a map $i: A \map B$ has the left lifting property
with respect to all special acyclic fibrations.
Use Lemma \ref{lem:factor1} to 
factor $i$ as a strict cofibration $i':A \map B'$ followed by
a special acyclic fibration $p:B' \map B$.  Then we have a square
\[
\xymatrix{
A \ar[r]^{i'} \ar[d]_{i} & B' \ar[d]^{p} \\
B \ar[r]_{=} & B,              }
\]
and a lift exists in this square by assumption.  Hence $i$ is a
retract of $i'$.  But retracts preserve strict cofibrations
by Proposition \ref{prop:retract}, so $i$ is again a strict cofibration.

Finally, suppose that $p: X \map Y$ has the right lifting property
with respect to all strict cofibrations.
Use Lemma \ref{lem:factor1} to 
factor $p$ as a strict cofibration $i:X \map X'$ followed by
a special acyclic fibration $p':X' \map Y$.  
Similarly to the previous paragraph, $p$ is a retract of $p'$.
\end{proof}

\begin{lemma}
\label{lem:lift2}
A map is an essentially levelwise acyclic cofibration if and only if 
it has the right lifting
property with respect to all strict fibrations.
Also, a map has the right lifting property with respect
to all essentially levelwise acyclic cofibrations if and only if 
it is a strict fibration.
\end{lemma}

\begin{proof}
The proof is the same as the proof of Lemma \ref{lem:lift1},
except that the roles of cofibrations and acyclic fibrations
are replaced by acyclic cofibrations and fibrations respectively.
Lemma \ref{lem:factor2} is relevant instead of Lemma \ref{lem:factor1}.
\end{proof}

\begin{proposition}
\label{prop:a-cofib}
A map is a strict acyclic cofibration if and only if it
is an essentially levelwise acyclic cofibration.
\end{proposition}

\begin{proof}
One direction follows from the definitions.

Let $i:A \map B$ 
be a strict weak equivalence and strict cofibration.  
We may assume that $i$ is a level representation that is a levelwise
weak equivalence.  Then we may use the argument of the proof of
Lemma \ref{lem:factor1} to factor $i$ as an essentially
levelwise acyclic cofibration $i':A \map B'$ followed by a
special acyclic fibration $p: B' \map B$.  Hence we have a square
\[
\xymatrix{
A \ar[r]^{i'} \ar[d]_i & B' \ar[d]^{p} \\
B \ar[r]_{=} & B,             }
\]
and a lift exists in this square by Lemma \ref{lem:lift1}
because $i$ is a strict cofibration and
$p$ is a special acyclic fibration.  Therefore, $i$ is a retract
of $i'$.  But retracts preserve essentially levelwise acyclic cofibrations
by Proposition \ref{prop:retract}, so $i$ is also
an essentially levelwise acyclic cofibration.
\end{proof}

\begin{proposition}
\label{prop:a-fib}
A map is a strict acyclic fibration if and only if it
is a retract of a special acyclic fibration.
\end{proposition}

\begin{proof}
A special acyclic fibration is a 
strict acyclic fibration by Lemma \ref{lem:spec-a-fib}.  For the other
direction, let $p:X \map Y$ 
be a strict weak equivalence and strict fibration.  
We may assume that $p$ is a level representation that is a levelwise
weak equivalence.  Then we may use the argument of the proof of
Lemma \ref{lem:factor1} to factor $p$ as an essentially
levelwise acyclic cofibration $i:X \map X'$ followed by a
special acyclic fibration $p': X' \map Y$.  
Similarly to the proof of Proposition \ref{prop:a-cofib},
$p$ is a retract of $p'$.
\end{proof}

\begin{theorem} \label{thm:strict}
Let $\C$ be a proper model category.  Then the classes of
strict weak equivalences, strict cofibrations, and strict fibrations
define a proper model structure on $\pro \C$.  
\end{theorem}

\begin{proof}
Completeness and cocompleteness follows from \cite[Prop.~11.1]{prospace}.
The two-out-of-three axiom is proved in Lemmas \ref{lem:composition}
and \ref{lem:2/3}.  The retract axiom for strict cofibrations and
strict weak equivalences follows from Proposition \ref{prop:retract}, where
it is shown that any class of essentially levelwise maps is always
closed under retract.
The retract axiom for strict fibrations is true by definition.

Using Propositions \ref{prop:a-cofib} and \ref{prop:a-fib},
the factorization axiom is given in Lemmas \ref{lem:factor1} and
\ref{lem:factor2}.
Similarly, the lifting axiom follows from 
Lemmas \ref{lem:lift1} and \ref{lem:lift2}.

This finishes the proofs of the basic model structure axioms.  It remains
to consider properness.
We shall show that the strict model structure is right proper;
the proof of left properness is dual.

Let $p: X \map Y$ be a strict fibration, and let $f: Z \map Y$ be a 
strict weak equivalence.  By Lemma \ref{lem:level-fib},
we know that $p$ is an essentially levelwise fibration.

We do not have a level representation
\[
\xymatrix@1{
Z \ar[r]^{f} & Y & X \ar[l]_{p}  }
\]
because we cannot necessarily represent $f$ by a levelwise weak equivalence
and $p$ by a levelwise fibration with the same indexing category.  
Nevertheless,
we do have a level representation
\[
\xymatrix@1{
Z \ar[r]^{f} & W \ar[r]^{g} & Y & \ar[l]_{p} X  }
\]
in which $p$ is a levelwise fibration, $f$ is a levelwise weak equivalence,
and $g$ is a pro-isomorphism (but not necessarily a levelwise isomorphism).
This gives us a level representation
\[
\xymatrix{
Z \times_{Y} X \ar[d] \ar[r]^{f'} & W \times_{Y} X \ar[r]^-{g'} \ar[d]_{p'} & 
  X \ar[d]^{p} \\
Z \ar[r]_{f} & W \ar[r]_{g} & Y,    }
\]
where the pullbacks are computed levelwise.
Because $g$ and therefore $g'$ are 
isomorphisms, 
$f'$ is a level representation for $Z \times_Y X \map X$.
Hence it suffices to show that $f'$ is a levelwise weak equivalence.

Since pullbacks preserve fibrations in $\C$, we know that
$p'$ is a levelwise fibration.  
Now the map $f'$ is a levelwise
pullback of a weak equivalence along a fibration, so it is a levelwise
weak equivalence because $\C$ is right proper.  
\end{proof}

\begin{remark}
\label{rem:proper}
The proof of properness in \cite[Prop.~17.1]{prospace} is incorrect, 
but the techniques
described in the above proof can be used to fix it.
\end{remark}

\subsection{Simplicial Model Structures}

If $\C$ is a simplicial category, then $\pro \C$ is again a 
simplicial category.
For any two pro-objects $X$ and $Y$, the \mdfn{simplicial mapping space
$\Map(X,Y)$} is defined to be $\lim_s \colim_t \Map_{\C} (X_t, Y_s)$.

Beware that the definitions of tensor and cotensor are 
straightforward for finite simplicial sets but are slightly subtle
in general.  If $X$ is a pro-object and $K$ is a finite simplicial
set, then $X \otimes K$ is defined by tensoring levelwise.  Similarly,
$X^K$ is defined by cotensoring levelwise.  

For an arbitrary simplicial set $K$, write it as a colimit 
$\colim_s K_s$, where each $K_s$ is a finite simplicial set.
Then define $X \otimes K$ to be $\colim^{\pros}_s (X \otimes K_s)$, where
the colimit is computed in the category $\pro \C$.
This is not the same as tensoring levelwise with $K$. 
Similarly, $X^K$ is defined to be $\lim^{\pros}_s (X^{K_s})$, where the
limit is computed in the category $\pro \C$.  Again, this
is not the same as cotensoring levelwise with $K$.

\begin{theorem}
\label{thm:simplicial}
If $\C$ is a proper simplicial model category, then
the strict structure on $\pro \C$ is also simplicial.
\end{theorem}

\begin{proof}
Most of the axioms are obvious or follow formally from the simplicial
structure on $\C$.  See \cite[\S 16]{prospace} for more details.
We shall show that if $i:K \map L$ is a cofibration of finite simplicial
sets and $j:A \map B$ is a strict cofibration in $\pro \C$, then
\[
A \otimes L \amalg_{A \otimes K} B \otimes K \map B \otimes L
\]
is a strict cofibration, and it is acyclic if either $i$ or $j$ is.
By adjointness \cite[Lem.~10.3.6]{PH}, 
this is equivalent to the usual SM7 axiom for
simplicial model categories.

We may assume that $i$ is a levelwise cofibration.  If $i$ is also acyclic,
then we may assume by Proposition \ref{prop:a-cofib} that $i$ is 
a levelwise acyclic cofibration.  Because $\C$
is a simplicial model category, the map
\[
A_s \otimes L \amalg_{A_s \otimes K} B_s \otimes K \map B_s \otimes L
\]
is a cofibration, and it is acyclic if either $i_s$ or $j$ is acyclic.
Because $K$ and $L$ are finite, this shows that
\[
A \otimes L \amalg_{A \otimes K} B \otimes K \map B \otimes L
\]
is a levelwise cofibration and that it is a levelwise acyclic cofibration
if either $i$ or $j$ is acyclic.
\end{proof}

\section{Non-Fibrantly Generated Model Categories}

In \cite[Cor.~19.3]{prospace}, it was shown that strict model structures
are not always cofibrantly generated \cite[Defn.~13.2.2]{PH}, even if
$\C$ is cofibrantly generated.  
For example, the strict model structure for pro-simplicial sets is
not cofibrantly generated.
In this section, we study whether
strict model structures are fibrantly generated.

Let $\lambda$ be any ordinal.  Then $\lambda$ is the partially ordered
set of all ordinals strictly less than $\lambda$.  
A \mdfn{$\lambda$-tower} $Z$ in a category is a contravariant functor
from $\lambda$ 
such that for all limit ordinals $\beta$, the object $Z_{\beta}$ is
isomorphic to $\lim_{\alpha < \beta} Z_{\alpha}$.
In other words, it is a diagram
\[
\cdots \map Z_{\beta} \map \cdots \map Z_1 \map Z_0
\]
of length $\lambda$.  Note that $Z_{\beta}$ is defined only for 
$\beta < \lambda$, not for $\beta = \lambda$.
This definition is dual to the notion of 
$\lambda$-sequence that arises in discussions of the small
object argument \cite[Defn.~12.2.1]{PH}.

\begin{definition}
\label{defn:cocell}
Let $C$ be a class of maps in a category.  A map
$p:X \map Y$ is a \mdfn{$C$-cocell complex} if there exists a
$\lambda$-tower $Z$ such that each map $Z_{\beta + 1} \map Z_{\beta}$
is a base change of a map belonging to $C$ and such that the
projection $\lim_{\beta < \lambda} Z_{\beta} \map Z_0$ is
isomorphic to the map $p$.
\end{definition}

This definition is dual to the usual notion of cell complex 
\cite[Defn.~12.4.7]{PH}.

\begin{proposition}
\label{prop:cocell}
Let $C$ be any class of maps in a category $\C$, and let $cC$ be the
image in $\pro \C$ of the class $C$ under the constant functor.  
If $f$ is a cofinite directed level representation such that each map 
$M_s f$ belongs to $C$, then $f$ is a $cC$-cocell complex.
\end{proposition}

\begin{proof}
Let $f$ be a cofinite directed level representation indexed by 
a cofinite directed set $I$ such that each $M_s f$ belongs to $C$.
Choose a well-ordering $\phi$ of $I$ that respects the ordering on $I$;
thus $\phi$ is a set isomorphism from $I$ to some ordinal $\lambda$ such
that $\phi(s) \geq \phi(t)$ whenever $s \geq t$.

We define a $\lambda$-tower $Z$ as follows.  Set $Z_0$ equal to $Y$.
For every successor ordinal $\beta + 1$, let $s$ be the unique element of
$I$ such that $\phi(s) = \beta$.  Define $Z_{\beta + 1}$ by the 
pullback square
\[
\xymatrix{
Z_{\beta + 1} \ar[r] \ar[d] & Z_{\beta} \ar[d] \\
cX_s \ar[r] & cY_s \times_{\lim^{\pros}_{t<s} cY_t} \lim^{\pros}_{t<s} cX_t. }
\]
Note that the bottom horizontal map is equal to $cM_s f$ because
$c$ commutes with finite limits.
Also note that by transfinite induction, each $Z_\beta$ comes
equipped with a map
\[
Z_\beta \map cY_s \times_{\lim^{\pros}_{t<s} cY_t} 
   \lim\nolimits^{\pros}_{t<s} cX_t.
\]

For limit ordinals $\beta$, define $Z_\beta$ to be 
$\lim^{\pros}_{\alpha < \beta} Z_\alpha$.
Thus we have constructed a $cC$-cocell complex $Z$.  A consideration
of universal properties shows that $\lim^{\pros}_{\beta < \lambda} Z_\beta$ 
is isomorphic to $X$; here we use that $X$ is isomorphic to
$\lim^{\pros}_s cX_s$.
Also, the projection
$\lim^{\pros}_{\beta < \lambda} Z_\beta \map Z_0$ is equal to $f$.
\end{proof}

\begin{proposition}
\label{prop:generate}
Let $\C$ be a proper model category, and let $I$ and $J$ be the classes of
fibrations and acyclic fibrations in $\C$ respectively.  
The class of strict fibrations is equal to the
class of retracts of $cI$-cocell complexes, and the 
class of strict acyclic fibrations is equal to the
class of retracts of $cJ$-cocell complexes.
\end{proposition}

\begin{proof}
We prove only the first statement; the proof of the second 
uses Proposition \ref{prop:a-fib} but is otherwise identical.

One implication follows immediately from Proposition \ref{prop:cocell}.
For the other direction, first 
observe that maps in $cI$ are special fibrations.
Base changes and compositions along $\lambda$-towers preserve right
lifting properties, so all $cI$-cocell complexes are strict fibrations.
Finally, retracts preserve lifting properties, so retracts of 
$cI$-cocell complexes are strict fibrations.
\end{proof}

\begin{remark}
One direction of the previous proof uses the
lifting property characterization of
strict fibrations, which we know because of Theorem \ref{thm:strict}, but
the other direction does not.
It seems plausible that for any category $\C$ and any class $C$ of
maps, a map in $\pro \C$ should be a retract of a $cC$-cocell complex
if and only if it is a retract of a map that has a cofinite directed
level representation for which all the relative matching maps belong to
$C$.  However, we have not been able to prove this claim in such generality.
\end{remark}

\begin{proposition} \label{prop:generate2}
Let $\C$ be a proper model category, and let $I$ and $J$ be the classes of
fibrations and acyclic fibrations in $\C$ respectively.  
The class of strict cofibrations is determined by the left 
lifting property with respect to $cJ$, and
the class of strict acyclic cofibrations is determined by the left 
lifting property with respect to $cI$.  
\end{proposition}

\begin{proof}
We only prove the first claim; the proof of the second is identical.

It is straightforward to check that every strict cofibration
has the left lifting property with respect to every element of $cJ$.
Conversely, if a map $i$ has the left lifting property with respect
to every element of $cJ$, then $i$ has the left lifting property
with respect to retracts of $cJ$-cocell complexes.
By Proposition \ref{prop:generate}, $i$ lifts with respect to all
strict acyclic fibrations, so it is a strict cofibration.
\end{proof}

Every object of every pro-category is cosmall \cite{duality}.
Therefore, Proposition \ref{prop:generate2} almost shows that
the strict model structure on $\pro \C$ is fibrantly generated.
The problem is that the collections $cI$ and $cJ$ of generating
strict fibrations and generating strict acyclic fibrations are
{\em not} sets.

\begin{proposition}
\label{prop:not-fib-gen}
The strict structure
on pro-simplicial sets is not fibrantly generated.  
\end{proposition}

\begin{proof}
Suppose that the strict structure on pro-simplicial sets is 
fibrantly generated.
Let $I$ be a set of generating strict fibrations.  By the
dual to the usual small object argument, 
every strict fibration is a retract of an $I$-cocell complex.

Now let $\lim I$ be the image of $I$ under the limit functor.
Since the constant functor $c$ and the limit functor form a Quillen
pair between $\C$ and $\pro \C$, 
every element of $\lim I$
is a fibration of simplicial sets.  

Every fibration
of simplicial sets is the image of a strict fibration under the limit functor,
so every fibration is a retract of a $(\lim I)$-cocell complex.
This means that the set $\lim I$ of fibrations detects acyclic cofibrations
of simplicial sets.  No such set exists, so we have obtained a
contradiction.
\end{proof}

\begin{remark}
\label{rem:not-fib-gen}
The fact that no set of fibrations detects
acyclic cofibrations of simplicial sets seems obvious, but the proof
is not elementary.  Bill Dwyer has shown us a proof, but it is too
lengthy to reproduce here.
\end{remark}

We describe 
one way to obtain a fibrantly generated model structure for pro-simplicial
sets.  Choose an uncountable cardinal $\kappa$.
We say that a simplicial set is $\kappa$-bounded if it has fewer than
$\kappa$ simplices.  

The category of $\kappa$-bounded simplicial sets equipped with the 
usual notions of weak equivalence, cofibration, and fibration is
almost a model category but not quite.  The problem is that not all
small limits and colimits exist.  Only the factorization axiom has a
non-obvious proof.  One must check that the usual small object
argument factorizations produce $\kappa$-bounded simplicial sets.

Now the category of $\kappa$-bounded simplicial sets is small and has
all finite limits.  This means that the category pro-($\kappa$-bounded
simplicial sets) has all small limits and colimits 
\cite[App.~4.3 and App.~4.4]{AM}.
We obtain a strict model structure
for pro-($\kappa$-bounded simplicial sets), and it is fibrantly generated.
As shown in Proposition \ref{prop:generate2}, the constant fibrations of
$\kappa$-bounded simplicial sets form a set of generating fibrations,
and the constant acyclic fibrations of 
$\kappa$-bounded simplicial sets form a set of generating acyclic fibrations.

Therefore, the strict model structure on pro-($\kappa$-bounded simplicial
sets) is fibrantly generated.  In some applications of pro-simplicial sets,
it suffices to choose a $\kappa$ larger than any of the simplicial
sets occurring in the application.  Thus, it is possible sometimes to
use a fibrantly generated model structure to study the
homotopy theory of pro-simplicial sets.

\end{document}